\pdfoutput=1
\RequirePackage{ifpdf}
\ifpdf 
\documentclass[pdftex]{sigma}
\else
\documentclass{sigma}
\fi

\begin{document}

\allowdisplaybreaks

\renewcommand{\thefootnote}{$\star$}

\renewcommand{\PaperNumber}{014}

\FirstPageHeading

\ShortArticleName{Courant Algebroids. A Short History}

\ArticleName{Courant Algebroids. A Short History\footnote{This
paper is a~contribution to the Special Issue ``Symmetries of Dif\/ferential Equations: Frames,
Invariants and~Applications''.
The full collection is available
at
\href{http://www.emis.de/journals/SIGMA/SDE2012.html}{http://www.emis.de/journals/SIGMA/SDE2012.html}}}

\Author{Yvette KOSMANN-SCHWARZBACH}

\AuthorNameForHeading{Y.~Kosmann-Schwarzbach}

\Address{Centre de Math\'ematiques Laurent Schwartz,
\'Ecole Polytechnique, F-91128 Palaiseau, France}
\Email{\href{mailto:yks@math.polytechnique.fr}{yks@math.polytechnique.fr}}
\URLaddress{\url{http://www.math.polytechnique.fr/cmat/kosmann/}}

\ArticleDates{Received December 03, 2012, in f\/inal form February 14, 2013; Published online February 19, 2013}

\Abstract{The search for a geometric interpretation of the constrained brackets
of Dirac led to the def\/inition of the Courant bracket. The search for
the right notion of a~``double'' for Lie bialgebroids led to the
def\/inition of Courant algebroids. We recount the emergence of these
concepts.}


\Classification{01A65; 53D17; 17B62; 17B66; 22A22; 53C15}

\begin{flushright}
\begin{minipage}{100mm}\it I dedicate this little memoir
to Peter Olver on the occasion\\ of his 60th birthday, with
friendship and in admiration.
\end{minipage}
\end{flushright}

\medskip

\renewcommand{\thefootnote}{\arabic{footnote}}
\setcounter{footnote}{0}

In  1986, in the historical notes to his {\it Applications of Lie Groups to
  Differential Equations} \cite{O},
Peter Olver pointed out that the concept of a ``Poisson structure''
 was already known to Sophus Lie, under the name ``function
 group,'' a century before
Andr\'e Lichnerowicz formally def\/ined Poisson manifolds
in 1977 \cite{Lichnerowicz}\footnote{See also Alan
  Weinstein's earlier historical note in {\it Expositiones
    mathematicae} \cite{AW}. For the development
of Poisson geometry up to 1998, see Weinstein's survey
\cite{AWsurvey}.
For more information on the history of Poisson brackets and Poisson
  geometry, see several chapters in the forthcoming book~\cite{add1}.}. Ten years before Lichnerowicz's
article, and in a completely unrelated inquiry,
Jean Pradines had def\/ined the new concept of Lie algebroid
as the inf\/initesimal counterpart of
Ehresmann's dif\/ferentiable groupoids, which are to-day called Lie groupoids
\cite{Pradines}. Both the geometry of Lichnerowicz's Poisson
manifolds
and the Ehresmann theory of groupoids
developed separately until the existence of a relationship
between the two theories was revealed when several mathematicians
and mathematical physicists independently
``discovered''
the Lie bracket of dif\/ferential $1$-forms on a Poisson
manifold\footnote{For a modern treatise on Lie groupoid and Lie
  algebroid theory, see \cite{Mack}. For elements of the history of the
  bracket of 1-forms on Poisson manifolds, see, e.g.,
footnote 5 in my survey~\cite{add2}.}.
In 1987, when Alan Weinstein, with Alain Coste and Pierre Dazord,
proved that the Lie algebroid of a symplectic groupoid is the
cotangent bundle of the base manifold equipped with this bracket of 1-forms
\cite{CDW}, it became clear that, for any Poisson manifold,
there is a Lie algebroid
structure on the cotangent bundle.
In fact, there is one more important property. The tangent and cotangent
bundles together constitute a Lie bialgebroid.
In 1988, Weinstein introduced the concept of a Poisson groupoid~\cite{AWPG}, and, in a joint paper published in 1994,
Kirill Mackenzie and Ping Xu def\/ined the concept of a Lie
bialgebroid as the
inf\/initesimal counterpart of Poisson groupoids \cite{MX}. A year later,
Kosmann-Schwarzbach proved that a pair of Lie algebroids in duality
form a Lie bialgebroid if and only if they def\/ine a pair of
dif\/ferential Gerstenhaber algebras\footnote{Gerstenhaber algebras are
  also called Schouten algebras. In fact, they are an abstract version of
  the algebra of multivector f\/ields on a manifold, equipped with the
  Schouten--Nijenhuis bracket.}, i.e., the Lie algebroid
structure of one def\/ines a~dif\/ferential which is a~derivation of the
Gerstenhaber bracket of
the exterior algebra of the other Lie algebroid \cite{yks3}.
A Poisson structure on a manifold thus gives rise to the prototypical
example of such a Lie
bialgebroid, consisting of the pair of the tangent and cotangent bundles.

Lie bialgebroids were ``just'' a generalization of the Lie bialgebras
that are the inf\/initesimal counterpart of the Poisson Lie groups.
The question that therefore arose was to def\/ine the double of such an
object, just as Vladimir Drinfeld had def\/ined the double of a Lie
bialgebra. His def\/inition, when $({\mathfrak
    g},{\mathfrak g}^*)$ is a Lie bialgebra, of the unique Lie algebra structure
  on ${\mathfrak g} \oplus {\mathfrak g}^*$ which leaves the canonical
symmetric bilinear form invariant, and restricts to the given Lie algebra structures on $\mathfrak g$
and ${\mathfrak g}^*$,
is contained in his 3-page 1983 paper~\cite{D0}\footnote{Drinfeld does not use
the word ``double''. The
  def\/inition of a ``Manin triple'' is to be found in his
landmark ``Quantum
  groups'' article which appeared in 1987 in the proceedings of the
  International Congress of Mathematicians (Berkeley,~1986). The
  actual lecture was delivered on very short notice by Pierre Cartier
  because Drinfeld was not allowed to travel to the United States for
  the conference.}.
The problem turned out to be more complicated than it
seemed, and it eventually led to the def\/inition of Courant algebroids by
Zhang-Ju Liu, Alan Weinstein and Ping Xu in~1997~\cite{LWX}.
The question that was posed had, as a special case,
the problem~--
which had been raised long before in  Poisson geometry~-- of
def\/ining a bracket on the
sections of the Whitney sum of the tangent and cotangent bundles of a Poisson
manifold, combining the  Lie bracket of vector f\/ields and the
above-mentioned bracket of 1-forms\footnote{I remember discussing
  this problem with Franco Magri {\it circa} 1990, and not f\/inding
a satisfying answer.}.

There was an even more ``elementary'' problem in dif\/ferential geometry,
how to def\/ine a~bracket with reasonable properties on the direct sum of
the vector space of vector f\/ields and the vector space of dif\/ferential
$1$-forms which restricts to the Lie bracket of vector f\/ields and
vanishes on dif\/ferential forms.
It was Theodore (Ted) Courant's achievement to def\/ine such a bracket in his
thesis, published in~1990~\cite{Courant}.
Below I shall describe the motivation for Courant's work, coming from Dirac's
theory of constraints, and the work that Irene Ya.~Dorfman de\-ve\-lo\-ped
independently, in the context of the
Hamiltonian structures, i.e., the inf\/inite-dimensional Poisson
structures, arising in f\/ield theory, in the late 1980's
and until her untimely death\footnote{Irene Dorfman (born 1948), one
  of the leading experts on integrable systems, died
  in Moscow in 1994. For her life and works, see
the obituary by Oleg I.~Mokhov, Sergei P.~Novikov and
Andrei K.~Pogrebkov~\cite{MNP}.}.

Dmitry Roytenberg advanced
the general theory of Courant algebroids further by founding
it on supermanifold theory in his thesis, completed in 1999, but
which has remained unpublished~\cite{RoytenbergThesis}, and in
his subsequent work~\cite{Roytenberg4, Roytenberg1,Roytenberg2}.
He adopted the ``cotangent philosophy"
of Kirill Mackenzie according to which, in particular,
given a Lie algebroid $A$, the cotangent bundle $T^*(A) \approx
T^*(A^*)$ is a more fundamental object than the bundle $A\oplus
A^*$.  He also made use of the supergeometry
approach that he had learned from Theodore (Ted) Voronov.

This dif\/ferent, more categorical
line of inquiry concerning the concept of a ``double'' for
Lie bialgebroids was explained by Mackenzie in
several publications, starting in 1992
\cite{Mackenzie1992,Mackenzie2000,Mackenzie1998}
and expounded in his article, f\/irst a preprint on arXiv in~2006,
submitted to Crelle's Journal in 2008, but f\/inally published only in~2011~\cite{Mackenzie2011}. (Double Lie algebroids are
the inf\/initesimal counterparts of double Lie groupoids, and Mackenzie's theory
includes many constructions besides the doubles of Lie bialgebroids.)
The question of the double of
a Lie bialgebroid was also treated by Ted Voronov in~\cite{Voronov1},
using a ``super'' approach and the notion of graded manifolds. In~\cite{Voronov2}, he characterized double Lie algebroids
in terms of supergeometry and
he proved that ``Roytenberg's and Mackenzie's pictures give the
same notion of a double of a Lie bialgebroid.''

\looseness=-1
The rich theory of the ``generalized
geometry'' of $TM \oplus T^*M$ was
developed, f\/irst by
Nigel Hitchin in~2003~\cite{Hitchin}, then by Marco Gualtieri in
his thesis submitted in 2003, available on arXiv in 2004, but
published much later as \cite{Gualtieri}, and it has become a domain
of research in its own right\footnote{The vector bunlde $TM \oplus
  T^*M$ has been variously
called the ``big tangent  bundle,'' by Izu Vaisman, and the
``Pontryagin bundle,'' by Hiroaki Yoshimura and
Jerrold E.~Marsden in~\cite{YM},
``because of its fundamental role
in the [geometric interpretation of] Pontryagin's maximum principle''
in control theory.
They point out that the Whitney sum, $TQ \oplus T^*Q$, was f\/irst
investigated in Lagrangian mechanics
by Ray Skinner and Raymond Rusk~\cite{SR}, who themselves refer to
earlier work by Mark Gotay and James N.~Nester.}.

In what follows I shall try to describe how and when a
non-skew\-symmetric bracket on a~Courant algebroid was
determined, and why it was eventually preferred in many circumstances
to the skewsymmetric version,
and I shall explain why it
is now called the Courant--Dorfman bracket, or sometimes just the Dorfman
bracket.

Alan Weinstein's
f\/irst article on Dirac structures,
with Ted Courant\footnote{Ted Courant is a grandson of the
  mathematician Richard Courant (1888--1972).},
his doctoral student at the time, ``Beyond Poisson
structures'', appeared in 1988 in the
proceedings of the ``Journ\'ees \mbox{lyonnaises} de la Soci\'et\'e
Math\'ematique de France (26--30 mai 1986) d\'edi\'ees \`a
\mbox{A.~Lichnerowicz''~\cite{CW}}. Motivated by the work of Robert
Littlejohn\footnote{This motivation was recalled in an e-mail
message I received from Alan Weinstein (July~23, 2012).}
on the Hamiltonian theory of guiding center
motion~\cite{Littlejohn79,Littlejohn84,Littlejohn81},
their idea was to def\/ine Poisson brackets
on subalgebras of the algebra of smooth functions on a smooth
manifold, such as the constrained brackets of Dirac, or the brackets of
functions constant on the characteristic foliation of a degenerate
2-form. They succeeded in interpreting geometrically Dirac's brackets as they
appeared in his {\it Lectures on
Quantum Mechanics}~\cite{Dirac} by setting up
a framework unifying Poisson and presymplectic structures.
What generalizes both
bivector f\/ields and dif\/ferential 2-forms on a manifold~$M$,
or rather their graphs, from $T^*M$ to $TM$ for bivectors, and from~$TM$ to
$T^*M$ for forms, are subbundles of $TM \oplus T^*M$ that are
maximally isotropic
(with respect to the canonical symmetric, f\/iberwise bilinear form),
which they called ``Dirac bundles''\footnote{See
Yoshimura and Marsden~\cite{YM} for elements of the history of
the introduction of Dirac structures on vector spaces and on manifolds,
and of their use in the study of degenerate Lagrangian systems and in
control theory.}.

The novelty was in the introduction of the direct sum of the tangent
and cotangent bundles. The dif\/f\/iculty lay in the
def\/inition of a 3-tensor on the subbundle whose
vanishing was the desired integrability condition that would reduce to
$[\pi,\pi]=0$, i.e.,
the vanishing of the Schouten--Nijenhuis bracket of~$\pi$,
when the subbundle is the graph of a bivector~$\pi$,
and would reduce to
${\rm d}\omega =0$, i.e., $\omega$ is closed, when the subbundle is
the graph of a 2-form $\omega$.
They solved this problem by
introducing a trilinear map $T$ on sections of $TM \oplus T^*M$
(see~\cite[p.~44--45]{CW}).
In the current literature, ``integrable Dirac subbundles'' are often
simply called
``Dirac bundles'', and they are said to def\/ine a ``Dirac structure''
on the base manifold.

Two years later, an important, additional novelty appeared in
Ted Courant's thesis, ``Dirac manifolds'', which was published in the {\it
  Transactions of the American Mathematical Society} in~1990~\cite{Courant}. He succeeded in def\/ining a skewsymmetric bracket on
$TM \oplus T^*M$, later called ``the Courant bracket'', and he def\/ined
a Dirac subbundle to be ``integrable'' if its space of sections is
closed under this bracket, in which case it is a Lie algebroid.
For sections $(X,\xi)$ and $(Y,\eta)$ of
  $TM \oplus T^*M$,
the Courant bracket is
\[
[(X,\xi), (Y,\eta)] = ([X,Y], {\mathcal L}_X\eta -{\mathcal L}_Y\xi
- \tfrac{1}{2} \, {\rm d} (i_X\eta - i_Y\xi)).
\]
Courant immediately stated, ``In general, this is {\it not} a
Lie-algebra bracket'', and this fact prompted some of the later
developments described below.
The 3-tensor $T$ of the joint
paper~\cite{CW} could be expressed simply
in terms of this new skewsymmetric bracket,
and Courant proved that the integrability of a Dirac subbundle is
equivalent to the vanishing of $T$, thus showing that his def\/inition
of ``integrability'' was equivalent to the one given in~\cite{CW}.
Therefore~$T$
appeared as the defect in the Jacobi identity for the
skewsymmetric Courant bracket, but this was not yet stated explicitly.

There arose the general question of how to def\/ine
a suitable notion of the double of a Lie bialgebroid.
In 1997,
Zhang-Ju Liu, Alan Weinstein and
Ping Xu proposed a solution to
this problem in their article, ``Manin triples for Lie
bialgebroids''~\cite{LWX}, by
def\/ining a skewsymmetric bracket on the space of sections of $A \oplus
A^*$, where $(A,A^*)$ is a Lie bialgebroid,
that generalizes the Courant bracket of $TM \oplus T^*M$.
They further introduced
a general notion of ``Courant algebroids'' (with a non-degenerate,
symmetric, f\/iberwise bilinear form and a skewsymmetric
bracket),
abstracted from the structure of the
doubles of Lie bialgebroids.
They wrote in their introduction:
``We found that if the bracket on a Courant algebroid is modif\/ied by
the addition of a symmetric term, many of the anomalies for the
resulting asymmetric [i.e., non-skewsymmetric]
bracket become zero,'' and
they asked, ``What is the geometric meaning of such asymmetric brackets?''
In a remark~\cite[p.~554]{LWX}, they introduced a non-skewsymmetric
bracket that they called ``a twisted bracket'', and listed three
of its properties. It is
easy to show that its skew-symmetrization is indeed the Courant
bracket.
They wrote, ``It would be nice to interpret
equation (i) [the Jacobi identity with a non-zero right-hand-side] in
terms of this twisted bracket. The geometric meaning of this twisted
bracket remains a mystery to us.'' All that was missing was the
interpretation of this bracket as a Loday bracket satisfying the
Jacobi identity in Leibniz form and, in fact, \v{S}evera and Weinstein
wrote later,
``It was observed [in 1998] by Kosmann-Schwarzbach, Xu, and \v Severa
(all unpublished) that
the non-skewsymmetric version of the bracket satisf\/ied the Jacobi
identity written in Leibniz form'' (\cite[p.~146]{SW}, also see
\cite[p.~527]{KW}).
Jean-Louis Loday\footnote{Jean-Louis Loday (born 1946) died in an
 accident at sea of\/f the coast of Brittany
in June~2012.}
 had introduced the concept of what he called
``Leibniz algebras'' in 1993~\cite{Loday}, but they have since been
called ``Loday algebras''. They are vector spaces equipped with a
non-skewsymmetric version of the Lie
brackets of Lie algebras that satisf\/ies
the following form of the Jacobi identity,
\[
[x,[y,z]] = [[x,y],z] + [y,[x,z]],
\]
which states that, for each element $x$, the adjoint map $[x, \cdot]$
is a derivation of the bracket, recalling the Leibniz rule for
the derivative of the product of two functions, hence the name ``Leibniz
algebra''.

In her e-mail message to Alan Weinstein of September~18, 1998,
Kosmann-Schwarzbach demonstrated that the Courant bracket
on $TM \oplus T^*M$ is  a
``derived bracket''\footnote{I~def\/ined the concept of derived bracket
in 1995,
inspired by unpublished notes that Jean-Louis Koszul had sent me,
and I published the def\/inition and properties of derived brackets with
applications to Poisson geometry in
the {\it Annales de l'Institut Fourier} in 1996 \cite{yks1}.
This
concept was also known to Ted Voronov,
who later developed the theory of higher derived brackets and their
relation to $L_\infty$-algebras.}. This
proof was published only~much later in her ``Derived
brackets'' paper \cite{yks2}, where the non-skewsymmetric bracket is
obtained in a natural way
as the derived bracket of the commutator of endomorphims of
the space of dif\/ferential forms
by the de Rham dif\/ferential, where vectors act by interior multiplication and
1-forms act by exterior multiplication. Then
the Courant bracket appears as the skew-symmetrization of this derived
bracket.

In Kosmann-Schwarzbach's
notes on her conversation with Pavol \v Severa at IH\'ES  on December~21,
1998, it is stated that together they had
verif\/ied that the Dirac structures of Courant and Weinstein
coincide with those of Dorfman as def\/ined in
1987 in~\cite{Dorfman1}. Later, in her book \cite{Dorfman2}, Dorfman
described  Dirac structures and applied them to the theory of
integrable equations.
She wrote in the introduction, ``Objects
called Dirac structures were introduced by Dorfman \cite{Dorfman1} as
natural algebraic analogues of f\/inite-dimensional structures f\/irst
introduced by Courant and Weinstein'', and she referred to their
work. Her citation, in both her article and her book,
of the work of Courant and Weinstein
as a Berkeley preprint dated 1986
was probably generous, because her
own work on ``the algebraic framework'' for Dirac structures
was, in fact, independent of theirs. We remark that,
although she
had f\/inally visited the West to participate in the workshop
on ``The Geometry of Hamiltonian Systems'' in Berkeley in June~1989,
she had not seen their article in published form, in the 1988 volume of the
``S\'eminaire Sud-Rhodanien de G\'eom\'etrie''~\cite{CW}, before she
completed her book.

It was Dmitry Roytenberg, another of Alan Weinstein's doctoral
students, who made
further progress in the theory of Courant algebroids in his thesis~\cite{RoytenbergThesis}.
He introduced the
non-skewsymmetric bracket,
gave a new def\/inition of Courant algebroids, containing f\/ive axioms,
in terms of this bracket,
and he proved the equivalence of
the old and the new def\/initions. These axioms can also be found in
the f\/irst of  {\v S}evera's e-mail letters to Weinstein of 1998~\cite{S}.
In~2001, {\v S}evera and Weinstein
published this def\/inition of Courant algebroids in~\cite{SW},
and Roytenberg f\/inally
published it in 2002 in
his ``Graded symplectic supermanifolds and Courant algebroids'' \cite{Roytenberg1}.
In his thesis, Roytenberg introduced a new interpretation of the
double of a Lie bialgebroid, as a homological\footnote{A vector f\/ield
  $X$ on a supermanifold is homological if $[X,X] =0$. It endows the
  manifold with the structure of a ``$Q$-manifold''.}
Hamiltonian vector f\/ield on an even symplectic supermanifold, thus
extending the def\/inition of a Lie algebroid structure as a homological vector
f\/ield, due to Arkady Vaintrob \cite{Vaintrob}. Combining the
``cotangent philosophy''
of Mackenzie with what he called
``Kosmann-Schwarzbach's picture of a Lie bialgebra''\footnote{See footnote~\ref{footnotebig} below.}, he recovered the
non-skewsymmetric bracket on $A \oplus A^*$ as the restriction of a
derived bracket of the canonical symplectic structure on a~graded version of the cotangent bundle, $T^*(A)
\approx T^*(A^*)$.

In 2002, Kyousuke Uchino showed that three of the axioms and one def\/ining
condition implied the other two axioms \cite{U}. A year later,
Janusz Grabowski and Giuseppe Marmo proposed a~def\/inition of the more general
Courant--Jacobi algebroids that required only four axioms~\cite{GM}.
Then in 2005, Kosmann-Schwarzbach proved that three of
Roytenberg's axioms for Courant algebroids
imply the other two~\cite{yks4}. One of these three axioms
is the Jacobi identity in Leibniz form for the
non-skewsymmetric bracket. The other two are expressed, as in
Roytenberg's thesis,
in a form that clearly shows that Courant
algebroids are a vector-bundle version of
the Lie algebras with an invariant symmetric
bilinear form (``quadratic Lie algebras'')\footnote{The
proof that the Leibniz rule and the morphism property of the
anchor are a consequence of these three axioms followed the same lines as
that of a property of Lie algebroids
in our paper with Franco Magri~\cite{add3}. We remark that the redundancy of some of
the axioms of Courant algebroids is proved using
the assumption that the symmetric
bilinear form is non-degenerate. In the more general case considered in~\cite{Bressler} and~\cite{Roytenberg3}, where no
such assumption is made, additional axioms are needed.}.

For the case of a
Courant algebroid which is the double of a Lie bialgebroid or,
more ge\-nerally, of a quasi-Lie bialgebroid, or of a proto-bialgebroid,
$(A,A^*)$, the formula for the Courant--Dorfman bracket is a
straightforward generalization of the case of the generalized tangent
bundles, using the ``big bracket''\footnote{\label{footnotebig}The
``big bracket'' for vector spaces
can be found in the paper by Bertram Kostant and Shlomo Sternberg~\cite{add4}. It
was f\/irst applied to the theory of Lie
bialgebras by Pierre Lecomte and Claude Roger in 1990, and I later
made extensive use of it, in~\cite{add5}, and in subsequent publications.}
on the algebra of functions on $T^*[2]A[1]$ as def\/ined
by Roytenberg, f\/irst in \cite{RoytenbergThesis}, then in his publication~\cite{Roytenberg2} in~2002, and the
derived bracket formula.
In particular, for the non-skewsymmetric bracket
on $TM \oplus T^*M$, the expression
of the derived bracket,
\[
[X+\xi, Y+ \eta] =\{\{X+\xi,
\mu\}, Y+ \eta\ \}= [X,Y] + {\mathcal L}_X\eta - i_Y {\mathrm d}\xi,
\]
for $X, Y \in \Gamma(TM)$, $\xi, \eta \in \Gamma(T^*M)$,
where $\mu$ is the Lie bracket of
vector f\/ields seen as a function on the supermanifold $T^*[2]TM[1]$,
coincides with the expression considered by Dorfman, in the context of complexes
over Lie algebras, in \cite[p.~242]{Dorfman1}, and in her
Theorem~2.1 in \cite{Dorfman2},
\[
[X+\xi, Y+ \eta] = [X,Y] + i_X {\mathrm d}\eta - i_Y {\mathrm d}\xi +
{\mathrm d}\langle X,\eta\rangle,
\]
in order to characterize Dirac structures. Whence, by extension,
the name ``Dorfman bracket'' that is now given
to the non-skewsymmetric bracket on any Courant algebroid.

In the more general case of an arbitrary Courant algebroid,
it was Roytenberg who proved in 2002 that the non-skewsymmetric Courant
bracket is a derived
bracket~\cite{Roytenberg1}.
To this end, he extended his own work on Lie bialgebroids in his
thesis and in~\cite{Roytenberg2}, def\/ining
a graded Poisson bracket on the ``minimal
symplectic realization'' of the bundle, based on
a construction that had been
suggested by Alan Weinstein in the Spring of 1999 and later simplif\/ied
by Pavol {\v S}evera in an unpublished letter to Weinstein
\cite[no.~7]{S}\footnote{In 2012, David Li-Bland
and Eckhard Meinrenken proposed a new construction of the
relevant Poisson bracket. See Shlomo Sternberg's
lectures at Harvard in 2012 and at the Conference in memory of
Jean-Marie Souriau,
Aix-en-Provence, June 2012.}.

Not only do Dirac structures have many applications in mechanics,
as developed in
the work of
Jerrold Marsden (1942--2010) and many others, in the AKSZ sigma-models
as was f\/irst shown by Noriaki Ikeda~\cite{Ikeda} (see~\cite{Roytenberg4}),
and even in supergravity as in the articles of
Daniel Waldram and
his co-authors (see, e.g.,~\cite{Waldram}),
many purely theoretical developments have taken place since 2002.
Twisted Courant algebroids,
also called Courant algebroids with background, were introduced.
{\v S}evera def\/ined a cohomology class, now called ``the {\v
  S}evera class''. In 2007, Paul Bressler def\/ined a Pontryagin class
in the generalized setting of transitive Lie algebroids,
and he showed that it is
an obstruction to the existence of a ``Courant extension''; in addition,
he related the theory of Courant
algebroids to conformal f\/ield theory~\cite{Bressler}.
Then, in 2009, there appeared the article by Roytenberg where he
def\/ined and studied an algebraic analogue of
Courant algebroids which he called ``Courant--Dorfman
algebras''~\cite{Roytenberg3}. Their relationship to Courant algebroids is
analogous to that of Lie--Rinehart algebras to Lie
algebroids\footnote{``Lie--Rinehart algebras'' is the name given by
  Johannes Huebschmann in 1990 to the algebraic counterpart of Lie
  algebroids. They are also called Lie pseudo-algebras.}.
The 2012 thesis by David
Li-Bland~\cite{add6} not only constitutes
a new contribution to the theory of Courant
algebroids, but also contains a useful list of references.

Symmetries are a fundamental feature of mathematical and physical
theories. Lie's conti\-nuous groups, which are now called Lie groups,
and their inf\/initesimal counterparts, which are now called
Lie algebras, have become too restrictive a framework for geometry,
algebra and mathematical physics. Lie groupoids
and Lie algebroids of\/fer a more general one, and Courant
algebroids have become a new, necessary concept in this wider framework,
where Drinfeld's double of a Lie bialgebra did not have an obvious
analogue. Lie algebroids and Courant algebroids are, in a sense,
inf\/initesimal objects. Lie algebroids correspond to Lie groupoids. Lie
bialgebras correspond to Poisson--Lie groups.
To what do Courant algebroids correspond?
All I know is that
Jean-Louis Loday had proposed to call these unknown objects
``coquecigrues''
and that, although some advances
have been made recently, the search
for coquecigrues\footnote{Rabelais, {\it Gargantua}, I.49. Picrochole
is told by an old hag ``que son royaulme
  luy serait rendu \`a la venue des cocquecigrues''.} is still on.

\subsection*{Acknowledgements}

\looseness=-1
I wish to thank Alan Weinstein as well as Dmitry Roytenberg and Pavol
{\v S}evera for their very useful correspondence and discussions
concerning their work on Courant algebroids. Thanks are also due to
Ted Voronov and to the referees for their
constructive remarks that were incorporated in the present version.
This short history was begun as a memorandum written
for Shlomo Sternberg, whom I thank for giving me the incentive to delve into
some near-contemporary history.

\pdfbookmark[1]{References}{ref}
\LastPageEnding

\end{document}